\documentclass{amsart}
\usepackage[english,francais]{babel}
\usepackage[dvips]{graphicx}
\usepackage{indentfirst}
\usepackage{fancyhdr}
\usepackage{amsfonts}
\usepackage{graphicx,epsfig}
\usepackage{euscript}
\DeclareMathAlphabet\mathscr{U}{eus}{m}{n}

\newtheorem{thm}{Theorem}
\newtheorem{prop}{Proposition}

\begin{document}
\pagestyle{headings}
\pagestyle{fancy}
\fancyhf{}
\fancyhead[RO,LE]{\thepage}
\fancyhead[CE]{Martin CELLI}
\fancyhead[CO]{The central configurations of four masses \(x\), \(-x\), \(y\), \(-y\)}

\begin{center}
{\large \bf The central configurations of four masses \(x\), \(-x\), \(y\), \(-y\).}\\
\(\)\\
Martin CELLI.\\
{\it Scuola Normale Superiore, Classe di Scienze.
Piazza dei Cavalieri 7, 56126 PISA. ITALY.}\\
{\it e-mail: m.celli@sns.it}\\ 
\(\)
\end{center}

{\bf Abstract.}
The configuration of a homothetic motion in the \(N\)-body problem
is called a central configuration. In this paper, we prove that there are exactly three
planar non-collinear central configurations for
masses \(x\), \(-x\), \(y\), \(-y\) with \(x\ne y\)
(a parallelogram and two trapezoids)
and two planar non-collinear central configurations for
masses \(x\), \(-x\), \(x\), \(-x\) (two diamonds).
Except the case studied here,
the only known case where the four-body central configurations
with non-vanishing masses can be listed is the case with equal masses
([1],[2]), which requires the use of a symbolic computation program.
Thanks to a lemma used in the proof of our result, we also show
that a co-circular four-body central configuration has
non-vanishing total mass or vanishing multiplier.\\ 


\section{Introduction.}

We are interested in configurations with \(N\) punctual
bodies which interact through gravitation, with masses
\(m_1\), ..., \(m_N\), and positions \(\vec{r}_1\), ..., \(\vec{r}_N\).
The masses are positive or negative but do not vanish,
the positions belong to a Euclidean vector space. The motion of the bodies is given by Newton's equations:
\[\ddot{\vec{r}}_i
=\vec{\gamma}_i(\vec{r}_1, ..., \vec{r}_N)
=\sum _{j\in \{1, ..., N\} \setminus \{i\}}
m_j\frac{\vec{r}_j-\vec{r}_i}{||\vec{r}_j-\vec{r}_i||^3} \cdot \]
We have:
\[\vec{\gamma} _i
=\frac{1}{m_i}\frac{\partial U}{\partial \vec{r}_i},\]
where \(U\) is the Newtonian potential:
\[U(\vec{r}_1, ..., \vec{r}_N)
=\sum _{1\le i<j\le N}\frac{m_im_j}{||\vec{r}_j-\vec{r}_i||} \cdot\]
We study a special class of configurations, which are called {\sl central
configurations}. A configuration is said to be central with {\sl multiplier}
\(\xi \) if, and only if, for any \(i\), \(j\):
\[\vec{\gamma}_j(\vec{r}_1, ..., \vec{r}_N)
-\vec{\gamma}_i(\vec{r}_1, ..., \vec{r}_N)
=\xi (\vec{r}_j-\vec{r}_i)\cdot \]
A {\sl motion} is said to be {\sl homothetic} if, and only if,
there is a scalar \(\alpha \) which depends on time, such that we have, for any \(i\), \(j\):
\[\vec{r}_j(t)-\vec{r}_i(t)=\alpha (t)(\vec{r}_j(0)-\vec{r}_i(0))\cdot \]
A {\sl motion} is said to be {\sl homothetic with fixed center} if,
and only if, there is is scalar \(\alpha \) which depends on time
and a constant vector \(\vec{\Omega}\), such that we have, for any \(i\):
\[\vec{r}_i(t)-\vec{\Omega}=\alpha (t)(\vec{r}_i(0)-\vec{\Omega})\cdot \]
It can be shown that the motion associated with a configuration and vanishing
initial velocities is homothetic if, and only if, the configuration is central.
It can be shown that it is homothetic with fixed center if, and only if,
the configuration is central with multiplier \(\xi \ne 0\) ([5]).
This motivates the study of central configurations.\\

The computation of central configurations is a difficult problem as soon as
\(N\ge 4\). Their finiteness up to similarities is the subject of
Smale's sixth problem for the 21st century ([14]). It has only been proved
recently in the case of four bodies with positive masses,
and the proof requires a computer ([8]).
On the other hand, there exists a continuum of
five-body central configurations when we allow a negative mass ([11]).
A continuum of four-body central configurations was also found
for charged particles with electrostatic (not gravitational) interaction,
which have positive masses and positive or negative charges
([4]).
Four-body central configurations with positive masses were also studied
numerically ([13]). Thanks to symbolic computation,
they can be listed when one mass vanishes
and the three others are equal ([10], [9])
and when the four masses are equal ([1],[2]).\\

In this article, we enumerate the four-body non-collinear planar central configurations
with masses \(x\), \(-x\), \(y\), \(-y\), where \(x\) and \(y\)
are two non-vanishing real numbers. The proof does not require the use
of numerical or symbolic computation. We also prove that 
a co-circular four-body central configuration has
non-vanishing total mass or vanishing multiplier.

\section{The field generated by a gravitational or electric dipole is one-to-one.}

The gravitational (respectively electric) field \(\vec{\gamma}\)
which is generated by two bodies with masses \(-1\) and \(1\)
(respectively with charges \(1\) et \(-1\))
and positions \(\vec{r}_1=(-1,0)\) and \(\vec{r}_2=(1,0)\)
in the plane, is given by:
\[\vec{\gamma}(\vec{r})=-\frac{\vec{r}_1-\vec{r}}{||\vec{r}_1-\vec{r}||^3}
+\frac{\vec{r}_2-\vec{r}}{||\vec{r}_2-\vec{r}||^3}\cdot \]

\begin{prop}
\label{dipole}
For every \(\vec{r}\), \(\vec{r}'\) such as
\(\vec{\gamma}(\vec{r})=\vec{\gamma}(\vec{r}')\),
we have: \(\vec{r}=\vec{r}'\) or \(\vec{r}=-\vec{r}'\).
\end{prop}

{\bf Proof.} Setting \(\vec{r}=(u,v)\), we have:
\[\gamma _u(u,v)=\frac{1+u}{((1+u)^2+v^2)^{3/2}}+\frac{1-u}{((1-u)^2+v^2)^{3/2}},
\gamma _v(u,v)=\frac{v}{((1+u)^2+v^2)^{3/2}}-\frac{v}{((1-u)^2+v^2)^{3/2}}\cdot \]
\[\frac{\partial \gamma _u}{\partial u}(u,v)
=\frac{v^2-2(u+1)^2}{((u+1)^2+v^2)^{5/2}}-\frac{v^2-2(1-u)^2}{((1-u)^2+v^2)^{5/2}}\cdot \]
\[\frac{\partial \gamma _u}{\partial v}(u,v)=\frac{\partial \gamma _v}{\partial u}(u,v)
=-3v\left (
\frac{u+1}{((u+1)^2+v^2)^{5/2}}
+\frac{1-u}{((1-u)^2+v^2)^{5/2}}
\right )\cdot \]
\[\frac{\partial \gamma _v}{\partial v}(u,v)
=\frac{(u+1)^2-2v^2}{((u+1)^2+v^2)^{5/2}}
+\frac{2v^2-(1-u)^2}{((1-u)^2+v^2)^{5/2}}\cdot \]
We set:
\[A(u,v)=\frac{(u+1)^2}{((u+1)^2+v^2)^{5/2}}-\frac{(1-u)^2}{((1-u)^2+v^2)^{5/2}} \mbox{, }
B(u,v)=\frac{v^2}{((u+1)^2+v^2)^{5/2}}-\frac{v^2}{((1-u)^2+v^2)^{5/2}}\cdot \]
We have:
\[\frac{\partial \gamma _u}{\partial u}(u,v)
\frac{\partial \gamma _v}{\partial v}(u,v)
=(B(u,v)-2A(u,v))(A(u,v)-2B(u,v))\]
\[=-2(A(u,v)-B(u,v))^2+A(u,v)B(u,v)\cdot \]
Let us denote by \(J\) the determinant of the Jacobian matrix of \(\vec{\gamma}\). We have:
\[J(u,v)=-\left ( \frac{\partial \gamma _u}{\partial v}(u,v)\right )^2
+\frac{\partial \gamma _u}{\partial u}(u,v)
\frac{\partial \gamma _v}{\partial v}(u,v)\cdot \]
\[J(u,v)=-9\left ( 
\frac{(u+1)v}{((u+1)^2+v^2)^{5/2}}
+\frac{(1-u)v}{((1-u)^2+v^2)^{5/2}}\right )^2
-2\left (
\frac{(u+1)^2-v^2}{((u+1)^2+v^2)^{5/2}}
+\frac{-(1-u)^2+v^2}{((1-u)^2+v^2)^{5/2}}
\right )^2\]
\[+\left (
\frac{(u+1)^2}{((u+1)^2+v^2)^{5/2}}-\frac{(1-u)^2}{((1-u)^2+v^2)^{5/2}}
\right )
\left(
\frac{v^2}{((u+1)^2+v^2)^{5/2}}-\frac{v^2}{((1-u)^2+v^2)^{5/2}}
\right )\cdot
\]
\[J(u,v)=-8\left ( 
\frac{(u+1)v}{((u+1)^2+v^2)^{5/2}}
+\frac{(1-u)v}{((1-u)^2+v^2)^{5/2}}\right )^2
-2\left (
\frac{(u+1)^2-v^2}{((u+1)^2+v^2)^{5/2}}
+\frac{-(1-u)^2+v^2}{((1-u)^2+v^2)^{5/2}}
\right )^2\]
\[-\left ( 
\frac{(u+1)v}{((u+1)^2+v^2)^{5/2}}
+\frac{(1-u)v}{((1-u)^2+v^2)^{5/2}}\right )^2\]
\[+\left (
\frac{(u+1)^2}{((u+1)^2+v^2)^{5/2}}-\frac{(1-u)^2}{((1-u)^2+v^2)^{5/2}}
\right )
\left(
\frac{v^2}{((u+1)^2+v^2)^{5/2}}-\frac{v^2}{((1-u)^2+v^2)^{5/2}}
\right )\cdot \]
\[J(u,v)
=-8\left ( 
\frac{(u+1)v}{((u+1)^2+v^2)^{5/2}}
+\frac{(1-u)v}{((1-u)^2+v^2)^{5/2}}\right )^2
-2\left (
\frac{(u+1)^2-v^2}{((u+1)^2+v^2)^{5/2}}
+\frac{-(1-u)^2+v^2}{((1-u)^2+v^2)^{5/2}}
\right )^2\]
\[-\frac{2(u+1)(1-u)v^2+(u+1)^2v^2+(1-u)^2v^2}{((u+1)^2+v^2)^{5/2}((1-u)^2+v^2)^{5/2}}
\cdot \]
\[J(u,v)
=-8\left ( 
\frac{(u+1)v}{((u+1)^2+v^2)^{5/2}}
+\frac{(1-u)v}{((1-u)^2+v^2)^{5/2}}\right )^2
-2\left (
\frac{(u+1)^2-v^2}{((u+1)^2+v^2)^{5/2}}
+\frac{-(1-u)^2+v^2}{((1-u)^2+v^2)^{5/2}}
\right )^2\]
\[-\frac{4v^2}{((u+1)^2+v^2)^{5/2}((1-u)^2+v^2)^{5/2}}\le 0\cdot \]
We have just proved that \(\vec{\gamma}\) is a submersion
for \(\vec{r}\ne \vec{0}\).\\

We have, for every \(\vec{r}\): \(\vec{\gamma}(-\vec{r})=\vec{\gamma}(\vec{r})\).
So proposition \ref{dipole} is equivalent to saying that
for every \(\vec{R}\), the following equation:
$$\vec{\gamma}(\vec{r})=\vec{\gamma}(\vec{R}) \eqno(1)$$
with unknown variable \(\vec{r}\) has exactly one solution when \(\vec{R}=\vec{0}\),
and two solutions when \(\vec{R}\ne \vec{0}\). This assertion is trivial
for \(\vec{R}=\vec{0}\). So it remains to prove that the restriction of
\(\vec{\gamma}\) to \(\mathbb{R}^2\setminus \{\vec{r}_1, \vec{0}, \vec{r}_2\}\)
is a \(2\)-fold covering map.\\

For \(\vec{R}=(1/2,0)\), we easily check that equation \((1)\)
has two solutions: \((-1/2,0)\) and \((1/2,0)\).
So it is enough to prove the following result:
if equation \((1)\) has exactly two solutions \(\vec{R}\) and \(-\vec{R}\)
at a point \(\vec{R}\in \mathbb{R}^2\setminus \{\vec{r}_1, \vec{0}, \vec{r}_2\}\),
this is still true in a neighborhood of \(\vec{R}\).
Let us assume this to be false. Then
there exists a sequence \((\vec{R}_n)_{n\in \mathbb{N}}\)
converging towards \(\vec{R}\) such that the number of solutions of
equation \((1)\) associated with \(\vec{R}_n\) is different from \(2\).
As \(d\vec{\gamma}(\vec{R})\) is invertible, there exists a neighborhood
\(\mathscr{U}\) of \(\vec{R}\) such that
for \(\tilde{R}\in \mathscr{U}\), equation \((1)\)
associated with \(\tilde{R}\) has exactly one solution \(\tilde{R}\) in \(\mathscr{U}\).
Then for \(n\) large enough,
the solutions in \(\mathscr{U}\cup (-\mathscr{U})\)
of equation \((1)\) associated with \(\vec{R}_n\)
are \(\vec{R}_n\) and \(-\vec{R}_n\).
This equation must have a third solution
\(\vec{R}_n'\). If \((\vec{R}_n')_{n\in \mathbb{N}}\) were bounded, we could
extract a subsequence converging towards a solution \(\vec{R}'\) of
equation \((1)\) associated with \(\vec{R}\). The solution \(\vec{R}'\)
would not belong to the open sets \(\mathscr{U}\)
and \(-\mathscr{U}\) and would be different
from \(\vec{R}\) and \(-\vec{R}\), which is impossible.
So the sequence \((\vec{R}_n')_{n\in \mathbb{N}}\) is not bounded. So we can assume
(if necessary, we can consider a subsequence): \(||\vec{R}_n'||\to +\infty\).
But this entails:
\(||\vec{\gamma}(\vec{R}_n)||=||\vec{\gamma}(\vec{R}_n')||\to 0\). Hence:
\(\vec{\gamma}(\vec{R})=\vec{0}\). Now we can check that this equality is false
for every \(\vec{R}\). QED

\section{Planar non-collinear central configurations with
masses \(x\), \(-x\), \(y\), \(-y\) and \(\xi =0\).}

\begin{thm}
\label{multzero}
Let \(x\) be a real number.
There is no planar non-collinear central configuration
with vanishing multiplier for the masses \(x\), \(-x\), \(x\), \(-x\).
Let \(x\) and \(y\) be two distinct non-vanishing real numbers.
There exists exactly one planar non-collinear central configuration
(up to similarities) with vanishing multiplier
for the masses \(x\), \(-x\), \(y\), \(-y\).
It is a parallelogram. The bodies with masses \(x\) and \(-x\)
(respectively the bodies with masses \(y\) and \(-y\))
are at the endpoints of a same diagonal.
The bodies whose masses have the smaller absolute value
are at the endpoints of the larger diagonal.
The two parallel sides whose endpoints have masses with the same sign are the largest ones.
\end{thm}

{\bf Proof.} Let \(x\) and \(y\) two non-vanishing real numbers.
We can assume: \(x\), \(y>0\). If necessary, we can exchange body \(1\) and body \(2\),
or body \(3\) and body \(4\).
We can also assume: \(x\le y\). If necessary, we can exchange
simultaneously body \(1\) and body \(3\), body \(2\) and body \(4\).\\

For any central configuration with vanishing multiplier, we have:
\(\vec{\gamma}_3(\vec{r}_1, ..., \vec{r}_4)
=\vec{\gamma}_4(\vec{r}_1, ..., \vec{r}_4)\).
This is equivalent to:
\(\vec{\gamma}(\vec{r}_3)=\vec{\gamma}(\vec{r}_4)\),
where \(\vec{\gamma}\) is
the function defined in the previous section.
According to proposition \ref{dipole},
the line segments \([\vec{r}_1,\vec{r}_2]\) and \([\vec{r}_3,\vec{r}_4]\)
have the same midpoint, so the configuration is a parallelogram.
The bodies with masses \(x\) et \(-x\)
(respectively the bodies with masses \(y\) and \(-y\))
are at the endpoints of a same diagonal.\\

Moreover, we have: \(\vec{\gamma}_1=\vec{\gamma}_3\). This is equivalent to:
\[(x+y)\frac{\vec{r}_1-\vec{r}_3}{||\vec{r}_1-\vec{r}_3||^3}
+(x-y)\frac{\vec{r}_3-\vec{r}_2}{||\vec{r}_3-\vec{r}_2||^3}
=x\frac{\vec{r}_1-\vec{r}_2}{||\vec{r}_1-\vec{r}_2||^3}
+y\frac{\vec{r}_4-\vec{r}_3}{||\vec{r}_4-\vec{r}_3||^3}\]
\[=x\frac{(\vec{r}_1-\vec{r}_3)-(\vec{r}_2-\vec{r}_3)}
{||\vec{r}_1-\vec{r}_2||^3}
+y\frac{(\vec{r}_1-\vec{r}_3)+(\vec{r}_2-\vec{r}_3)}
{||\vec{r}_4-\vec{r}_3||^3}\cdot\]
As the bodies are assumed not to be collinear,
the vectors \(\vec{r}_1-\vec{r}_3\) and \(\vec{r}_2-\vec{r}_3\)
are not collinear. So this is equivalent to:
\[\left \{ \begin{array}{c} 
\frac{2x}{||\vec{r}_2-\vec{r}_1||^3}
=\frac{x+y}{||\vec{r}_3-\vec{r}_1||^3}
+\frac{x-y}{||\vec{r}_3-\vec{r}_2||^3}\\
\frac{2y}{||\vec{r}_4-\vec{r}_3||^3}
=\frac{x+y}{||\vec{r}_3-\vec{r}_1||^3}
+\frac{y-x}{||\vec{r}_3-\vec{r}_2||^3}
\end{array} \right .\]
Let us assume: \(x=y\). These identities provide:
\(||\vec{r}_3-\vec{r}_1||=||\vec{r}_2-\vec{r}_1||
=||\vec{r}_4-\vec{r}_3||\).
This is impossible for a non-flat parallelogram.
So there is no planar non-collinear central configuration
with vanishing multiplier for the masses \(x\), \(-x\), \(x\), \(-x\).\\

From now on we assume: \(x<y\). As we are looking for configurations up to similarities,
we can also assume: \(||\vec{r}_3-\vec{r}_1||^2+||\vec{r}_3-\vec{r}_2||^2=1\).
Thus we have: \(||\vec{r}_3-\vec{r}_1||\), \(||\vec{r}_3-\vec{r}_2||<1\),
and the previous system provides:
\(\frac{2y}{||\vec{r}_4-\vec{r}_3||^3}>2y\),
which is equivalent to: \(||\vec{r}_4-\vec{r}_3||<1\).
Let us write the parallelogram law:
\[||\vec{r}_2-\vec{r}_1||^2+||\vec{r}_4-\vec{r}_3||^2
=2(||\vec{r}_3-\vec{r}_1||^2+||\vec{r}_3-\vec{r}_2||^2)=2\cdot \]
As \(||\vec{r}_4-\vec{r}_3||<1\), we necessarily have:
\(||\vec{r}_4-\vec{r}_3||<1<||\vec{r}_2-\vec{r}_1||\).
So the bodies whose masses have the smaller absolute value
are at the endpoints of the larger diagonal.
From this inequality we deduce:
\[\frac{2xy}{||\vec{r}_2-\vec{r}_1||^3}<\frac{2xy}{||\vec{r}_4-\vec{r}_3||^3}\cdot \]
\[y\left (\frac{x+y}{||\vec{r}_3-\vec{r}_1||^3}
+\frac{x-y}{||\vec{r}_3-\vec{r}_2||^3}\right )
<x\left (\frac{x+y}{||\vec{r}_3-\vec{r}_1||^3}
+\frac{y-x}{||\vec{r}_3-\vec{r}_2||^3}\right )\cdot \]
Hence: \(||\vec{r}_3-\vec{r}_2||<||\vec{r}_3-\vec{r}_1||\).
So the two parallel sides whose endpoints have masses with the same sign are
the largest ones.\\

Let us set:
\[(u,v)=\left (
\frac{1}{||\vec{r}_3-\vec{r}_1||^3},
\frac{1}{||\vec{r}_3-\vec{r}_2||^3} \right ),
(u',v')=\left (
\frac{1}{||\vec{r}_2-\vec{r}_1||^3},
\frac{1}{||\vec{r}_4-\vec{r}_3||^3} \right )\cdot \]
With these notations, the previously established system has the following expression:
\[\left \{ \begin{array}{c}
u'=\frac{x+y}{2x}u+\frac{x-y}{2x}v\\
v'=\frac{x+y}{2y}u+\frac{y-x}{2y}v
\end{array} \right .\]
The condition:
\(||\vec{r}_3-\vec{r}_1||^2+||\vec{r}_3-\vec{r}_2||^2=1\)
has the expression: \(u^{-2/3}+v^{-2/3}=1\).
The parallelogram law
has the expression: \(u'^{-2/3}+v'^{-2/3}=2\), which is equivalent to:
\[\left (\frac{x+y}{2x}u+\frac{x-y}{2x}v \right )^{-2/3}
+\left (\frac{x+y}{2y}u+\frac{y-x}{2y}v \right )^{-2/3}=2\cdot \]
We are going to show that the system \((\mathscr{S})\) which is defined
by these two equations:
\[\left \{ \begin{array}{c} 
f(u,v)=u^{-2/3}+v^{-2/3}=1\\
g(u,v)=\left (\frac{x+y}{2x}u+\frac{x-y}{2x}v \right )^{-2/3}
+\left (\frac{x+y}{2y}u+\frac{y-x}{2y}v \right )^{-2/3}=2
\end{array} \right .\]
has exactly one solution \((u,v)\).\\

Let \(\Gamma \) be the curve with equation: \(g(u,v)=2\).
Let us set:
\[h(u,v)=\left (\frac{x+y}{2x}u+\frac{x-y}{2x}v,
\frac{x+y}{2y}u+\frac{y-x}{2y}v\right )\cdot \]
The set \(\Gamma \) is the preimage of the curve with equation: \(u^{-2/3}+v^{-2/3}=2\)
by the invertible linear function \(h\). This curve is diffeomorphic to a line, so
\(\Gamma \) is diffeomorphic to a line.\\

\begin{figure}[h]
  \centering
     \includegraphics[width=7cm]{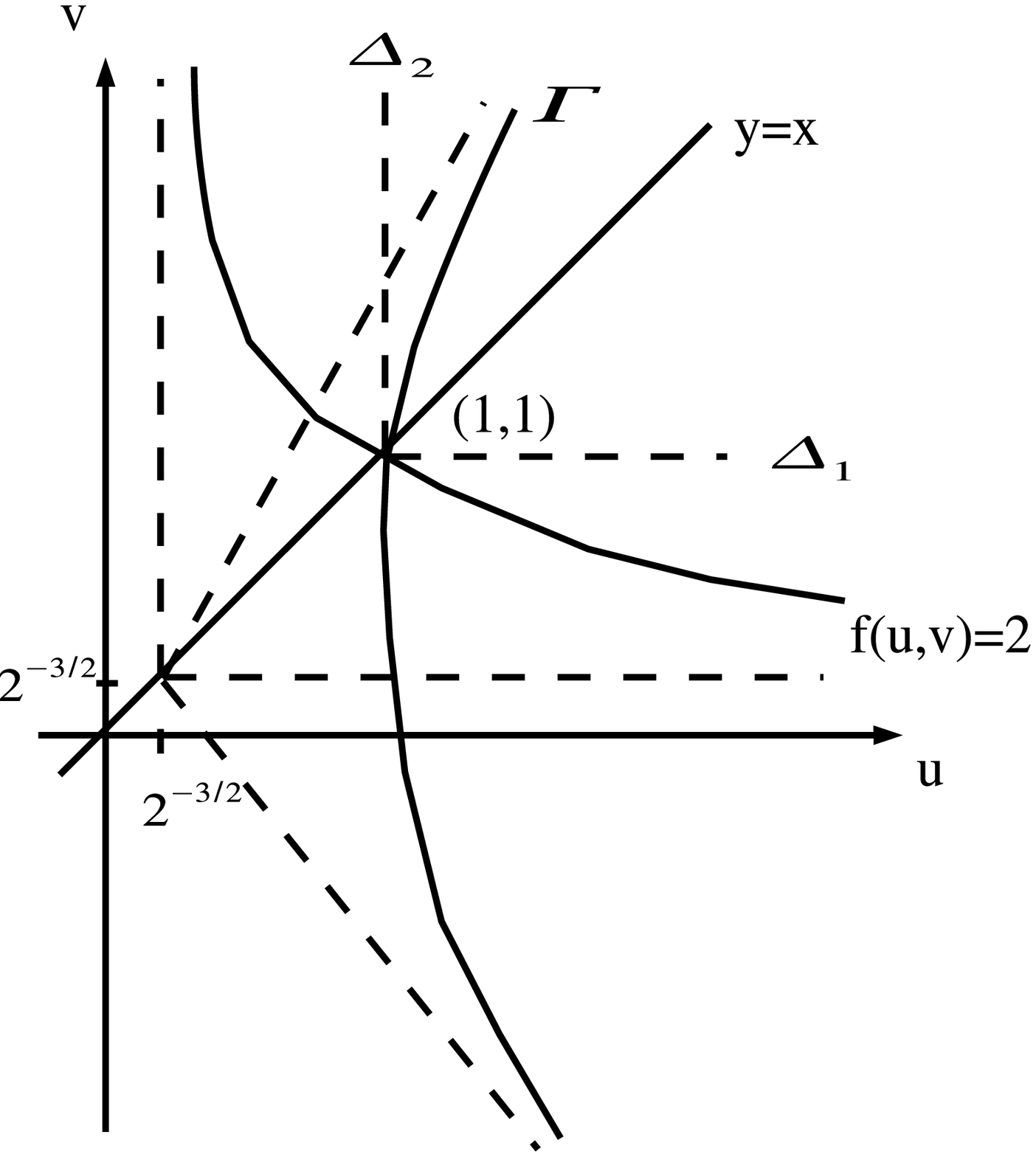}
       \caption{}
       \label{graphmult0}
       \end{figure}

Let \(\Delta _1\) be the ray defined by: \(v=1\), \(u\ge 1\).
Let \(\Delta _2\) be the ray defined by: \(u=1\), \(v\ge 1\).
Let us set: \(\Delta =\Delta _1 \cup \Delta _2\).
Let \(\gamma \) be the set of the \((u,v)\) such that:
\(u\), \(v>1\) and \(g(u,v)=2\).
Let us set: \(G_1(u)=g(u,1)\).
The function \(G_1\) is defined on \([1,+\infty[\)
and strictly decreasing. Now \(G_1(1)=2\),
so \(\Gamma \cap \Delta _1=\{(1,1)\}\).
Let us set: \(G_2(v)=g(1,v)\).
The function \(G_2\) is defined on \([1,(x+y)/(y-x)[\), and we have:
\[\frac{dG_2}{dv}(v)=\frac{y-x}{3}
\left (
\frac{1}{x}\left (\frac{x+y}{2x}+\frac{x-y}{2x}v \right )^{-5/3}
-\frac{1}{y}\left (\frac{x+y}{2y}+\frac{y-x}{2y}v \right )^{-5/3}
\right )\cdot \]
It can be checked that \(dG_2/dv>0\) if, and only if:
\[v>\frac{(x+y)(y^{2/5}-x^{2/5})}{(y-x)(x^{2/5}+y^{2/5})}\cdot \]
Now it easy to see that this inequality is always true
if \(y>x\) and \(v\ge 1\). So \(G_2\) is strictly increasing.
As \(G_2(1)=2\), we have:
\(\Gamma \cap \Delta _2=\{(1,1)\}\).
Thus: \(\Gamma \cap \Delta =\{ (1,1)\}\).
The set \(\Gamma \) is diffeomorphic to a line,
so \(\Gamma \setminus \{(1,1)\}\) has two connected components,
which are diffeomorphic to lines. As \(\Gamma \cap \Delta =\{ (1,1)\} \),
each connected component is included in a side of \(\Delta \).
We can check that only one of the two components is in the quarter plane defined by the inequalities: \(u>1\) et \(v>1\). So this connected component is \(\gamma \). Thus
\(\gamma \) is diffeomorphic to a line, an endpoint is \((1,1)\), the other is at infinity
(figure \ref{graphmult0}).\\

When \(u\) and \(v\) tend to infinity, \(f(u,v)\) tends to \(0<1\).
On the other hand, \(f(1,1)=2>1\). The system \((\mathscr{S})\)
has at least one solution.\\

The arguments previously used in order to prove the inequality:
\(||\vec{r}_4-\vec{r}_3||<1<||\vec{r}_2-\vec{r}_1||\),
enable to prove that:
\[v'=\frac{x+y}{2y}u+\frac{y-x}{2y}v>1>u'=\frac{x+y}{2x}u+\frac{x-y}{2x}v\]
on \(\gamma \). Hence:
\[\frac{1}{y}\left (\frac{x+y}{2y}u+\frac{y-x}{2y}v\right )^{-5/3}
<\frac{1}{x}\left (\frac{x+y}{2x}u+\frac{x-y}{2x}v\right )^{-5/3}\cdot \]
Thus on \(\gamma \) we have:
\[\frac{\partial g}{\partial v}(u,v)=
\frac{y-x}{3}\left (
\frac{1}{x}\left (\frac{x+y}{2x}u+\frac{x-y}{2x}v\right )^{-5/3}
-\frac{1}{y}\left (\frac{x+y}{2y}u+\frac{y-x}{2y}v\right )^{-5/3}
\right )>0\cdot \]
We immediately check: 
\[\frac{\partial f}{\partial u}(u,v)<0\mbox{, }
\frac{\partial f}{\partial v}(u,v)<0\mbox{, }
\frac{\partial g}{\partial u}(u,v)<0\cdot \]
So \(df\) and \(dg\) are collinear on \(\gamma \).
Thus the restriction of \(f\) to \(\gamma \) is strictly monotonic.
There exists exactly one point \((u,v)\) of \(\gamma \) such that
\(f(u,v)=1\). This point is the only solution of \((\mathscr{S})\),
which depends on \(x\) and \(y\) in a continuous way.
So there exists at most one planar non collinear central configuration
(up to similarities) with vanishing multiplier for the masses
\(x\), \(-x\), \(y\), \(-y\).\\

We are now going to show that the values of \(u\), \(v\), \(u'\), \(v'\) that we have just obtained actually define a non-collinear central configuration with vanishing multiplier.
Let us hold \(y\) constant. According to the following identity:
\[0<u'=\frac{x+y}{2x}u-\frac{y-x}{2x}v,\]
we have:
\[\frac{y-x}{x+y}v<u \cdot \]
Moreover: \(u<v\) (as
we established: \(||\vec{r}_3-\vec{r}_2||<||\vec{r}_3-\vec{r}_1||\)).
So \(u/v\to 1\) when \(x\to 0\).
According to the following identity:
\[\frac{v'}{v}=\frac{x+y}{2y}\frac{u}{v}+\frac{y-x}{2y},\]
we also have: \(v'/v\to 1\) when \(x\to 0\).
So the homothety class of \((u^{-1/3},v^{-1/3},v'^{-1/3})\)
tends to the class of the lengths of the sides of
an equilateral triangle when \(x\to 0\).
Thus, for \(x\) close to \(0\),
the numbers \((u^{-1/3},v^{-1/3},v'^{-1/3})\) are the lengths
\((||\vec{r}_3-\vec{r}_1||,||\vec{r}_4-\vec{r}_1||,||\vec{r}_4-\vec{r}_3||)\)
of a non-flat triangle, which is close to an equilateral triangle.\\

Let us now assume that for some \(x<y\), the numbers \((u^{-1/3},v^{-1/3},v'^{-1/3})\)
are not the lengths of the sides of a non-flat triangle.
Let us denote by \(u(t)\), \(v(t)\), \(v'(t)\) the values of \(u\), \(v\), \(v'\)
associated with \(x(t)=tx\) and \(y(t)=y\). For \(t\) close to \(0\), 
we have just seen that \(((u(t))^{-1/3},(v(t))^{-1/3},(v'(t))^{-1/3})\)
actually defines a non-flat triangle, whereas
\(((u(1))^{-1/3},(v(1))^{-1/3},(v'(1))^{-1/3})\)
does not define a non-flat triangle any more. 
The solutions \(u\), \(v\), \(v'\) depend on \(x\) and \(y\)
in a continuous way. So at a certain time \(t\le 1\), the triangle is flat.
Let us first assume that
\(\vec{r}_3\in ]\vec{r}_1,\vec{r}_4[\) at this time.
Let us set: \(a=(u(t))^{-1/3}=||\vec{r}_3-\vec{r}_1||\),
\(b=(v'(t))^{-1/3}=||\vec{r}_4-\vec{r}_3||\). We have:
\[||\vec{r}_4-\vec{r}_1||=(v(t))^{-1/3}=a+b,\]
\[(u'(t))^{-2/3}=g(u(t),v(t))-(v'(t))^{-2/3}=2-b^2=2f(u(t),v(t))-b^2\]
\[=2((u(t))^{-2/3}+(v(t))^{-2/3})-b^2=2(a^2+(a+b)^2)-b^2=(2a+b)^2\cdot \]
Hence:
\[\left \{ \begin{array}{c} 
\frac{2x(t)}{(2a+b)^3}
=\frac{x(t)+y(t)}{a^3}
+\frac{x(t)-y(t)}{(a+b)^3}\\
\frac{2y(t)}{b^3}
=\frac{x(t)+y(t)}{a^3}
+\frac{y(t)-x(t)}{(a+b)^3}
\end{array} \right .\]
Eliminating \(x(t)\) and \(y(t)\), we obtain:
\[\left ( \frac{1}{a^3}-\frac{1}{(a+b)^3} \right )^2
=\left ( \frac{1}{a^3}+\frac{1}{(a+b)^3}-\frac{2}{(2a+b)^3} \right )
\left ( \frac{1}{a^3}+\frac{1}{(a+b)^3}-\frac{2}{b^3} \right )\cdot \]
For sake of convenience, we can study this equation under the assumption: \(a=1\).
Then we obtain:
\[(b^3+(b+2)^3)((b+1)^3+1)=2(b^3(b+2)^3+(b+1)^3)\cdot \]
\[2b^6+12b^5+36b^4+66b^3+72b^2+48b+16
=2b^6+12b^5+24b^4+18b^3+6b^2+6b+2\cdot \]
No positive value of \(b\) is a solution. So the triangle cannot be flat.
We can obtain the same result for \(\vec{r}_4\in ]\vec{r}_1, \vec{r}_3[\).
Lastly, as we have:
\[v'(t)=\frac{x(t)+y(t)}{2y(t)}u(t)+\frac{y(t)-x(t)}{2y(t)}v(t),\]
we cannot simultaneously have: \(u(t)>v'(t)\) and \(v(t)>v'(t)\).
Hence: \(||\vec{r}_3-\vec{r}_1||\ge ||\vec{r}_4-\vec{r}_3||\)
or \(||\vec{r}_4-\vec{r}_1||\ge ||\vec{r}_4-\vec{r}_3||\).
So we cannot have: \(\vec{r}_1\in ]\vec{r}_3, \vec{r}_4[\).
Thus, the numbers \(((u(t))^{-1/3},(v(t))^{-1/3},(v'(t))^{-1/3})\)
cannot define a flat triangle. So the numbers
\((u^{-1/3},v^{-1/3},v'^{-1/3})=((u(1))^{-1/3},(v(1))^{-1/3},(v'(1))^{-1/3})\)
are the mutual distances
\((||\vec{r}_3-\vec{r}_1||,||\vec{r}_4-\vec{r}_1||,||\vec{r}_4-\vec{r}_3||)\)
of a non-flat triangle.\\

Let \(\vec{r}_2\) be the point such that
\((\vec{r}_1,\vec{r}_2,\vec{r}_3,\vec{r}_4)\)
is a parallelogram. We necessarily have:
\(\vec{\gamma}_3=\vec{\gamma}_4\). Moreover:
\[\left \{ \begin{array}{c} 
\frac{2x}{||\vec{r}_2-\vec{r}_1||^3}
=\frac{x+y}{||\vec{r}_3-\vec{r}_1||^3}
+\frac{x-y}{||\vec{r}_3-\vec{r}_2||^3}\\
\frac{2y}{||\vec{r}_4-\vec{r}_3||^3}
=\frac{x+y}{||\vec{r}_3-\vec{r}_1||^3}
+\frac{y-x}{||\vec{r}_3-\vec{r}_2||^3}
\end{array} \right .\]
This enables to prove that \(\vec{\gamma}_1=\vec{\gamma}_3\). Lastly:
\[m_2\vec{\gamma}_2
=-(m_1\vec{\gamma}_1+m_3\vec{\gamma}_3+m_4\vec{\gamma}_4)
=-(m_1+m_3+m_4)\vec{\gamma}_1=m_2\vec{\gamma}_1\cdot \]
As \(\vec{r}_1\), \(\vec{r}_3\) and \(\vec{r}_4\) are not collinear,
the parallelogram is not flat. So the configuration is non-collinear, central with vanishing multiplier for the masses \(x\), \(-x\), \(y\), \(-y\). QED

\section{Planar non-collinear central configurations with
masses \(x\), \(-x\), \(y\), \(-y\) and \(\xi \ne 0\).}

Let us define the vector of inertia
\(\vec{\lambda}\) of a \(N\)-body system with vanishing total mass (\(M=m_1+...+m_N=0\)) by:
\[\vec{\lambda}=m_1(\vec{r}_1-\vec{a})+...+m_N(\vec{r}_N-\vec{a})\]
(the identity does not depend on the origin \(\vec{a}\)).

\begin{prop}
\label{inertiestat}
For any central configuration with vanishing total mass:
\(\xi \vec{\lambda}=\vec{0}\).
\end{prop}

{\bf Proof.} We have to write:
\[\vec{0}=\sum _{i=1} ^N m_i\vec{\gamma}_i(\vec{r}_1, ..., \vec{r}_N)
=\sum _{i=1} ^N \xi m_i(\vec{r}_i-\vec{a})
=\xi \vec{\lambda}\cdot\]
QED\\

For every \(N\)-body configuration with dimension \(N-2\), there exists exactly one
\(N\)-tuple up to homotheties \((\Delta _1,...,\Delta _N)\)
such that \(\Delta _1+...+\Delta _N=0\)
and \(\Delta _1(\vec{r}_1-\vec{a})+...+\Delta _N(\vec{r}_N-\vec{a})=\vec{0}\)
(the identity does not depend on the origin \(\vec{a}\)).
For every configuration with dimension \(N-2\)
and vanishing \(\vec{\lambda}\), the \(N\)-tuples
\((m_1, ..., m_N)\) and \((\Delta _1, ..., \Delta _N)\) are collinear.

\begin{prop}
\label{cc4lambda0}
A planar non-collinear four-body configuration is central with vanishing
vector of inertia for a system of non-vanishing masses whose sum vanishes
if, and only if:\\
- No three-body subconfiguration is collinear.\\
- The mutual distances satisfy the following identity,
which does not involve the masses:
\[\frac{1}{||\vec{r}_2-\vec{r}_1||^3}+\frac{1}{||\vec{r}_4-\vec{r}_3||^3}
=\frac{1}{||\vec{r}_3-\vec{r}_1||^3}+\frac{1}{||\vec{r}_4-\vec{r}_2||^3}
=\frac{1}{||\vec{r}_3-\vec{r}_2||^3}+\frac{1}{||\vec{r}_4-\vec{r}_1||^3}\cdot \]
Then \((m_1,m_2,m_3,m_4)\) and \((\Delta _1,\Delta _2,\Delta _3,\Delta _4)\)
are collinear.
\end{prop}

{\bf Proof.} For a non-collinear four-body configuration with vanishing \(\vec{\lambda}\),
the \(N\)-tuples \((m_1,m_2,m_3,m_4)\)
and \((\Delta _1,\Delta _2,\Delta _3,\Delta _4)\) are collinear.
So \(\Delta _i\ne 0\) for every \(i\), and no three-body subconfiguration is collinear.
The configuration is central if, and only if, it satisfies
the Laura-Andoyer equations ([3]).
Thanks to the collinearity of \((m_1,m_2,m_3,m_4)\)
and \((\Delta _1,\Delta _2,\Delta _3,\Delta _4)\),
these equations are equivalent to:
\[\frac{1}{||\vec{r}_2-\vec{r}_1||^3}+\frac{1}{||\vec{r}_4-\vec{r}_3||^3}
=\frac{1}{||\vec{r}_3-\vec{r}_1||^3}+\frac{1}{||\vec{r}_4-\vec{r}_2||^3}
=\frac{1}{||\vec{r}_3-\vec{r}_2||^3}+\frac{1}{||\vec{r}_4-\vec{r}_1||^3}\cdot \]
QED\\

We can explain the degeneracy of the Laura-Andoyer equations
by noticing that the condition: \(\vec{\lambda}=\vec{0}\) enables to express the position of
a body as a barycenter of the three others.
Thus we have to deal with a three-body problem, with a different expression for the accelerations. Let us notice that this property is only the static variant of a more general
fact: for \(M=0\), the vector \(\vec{\lambda}\) is a first integral which
replaces the center of inertia. As it is invariant under translations,
the \(N\)-body problem becomes 'more integrable' ([5],[6]).\\

\begin{figure}[h]
  \centering
    \hspace{-4cm}
     \includegraphics[width=10cm]{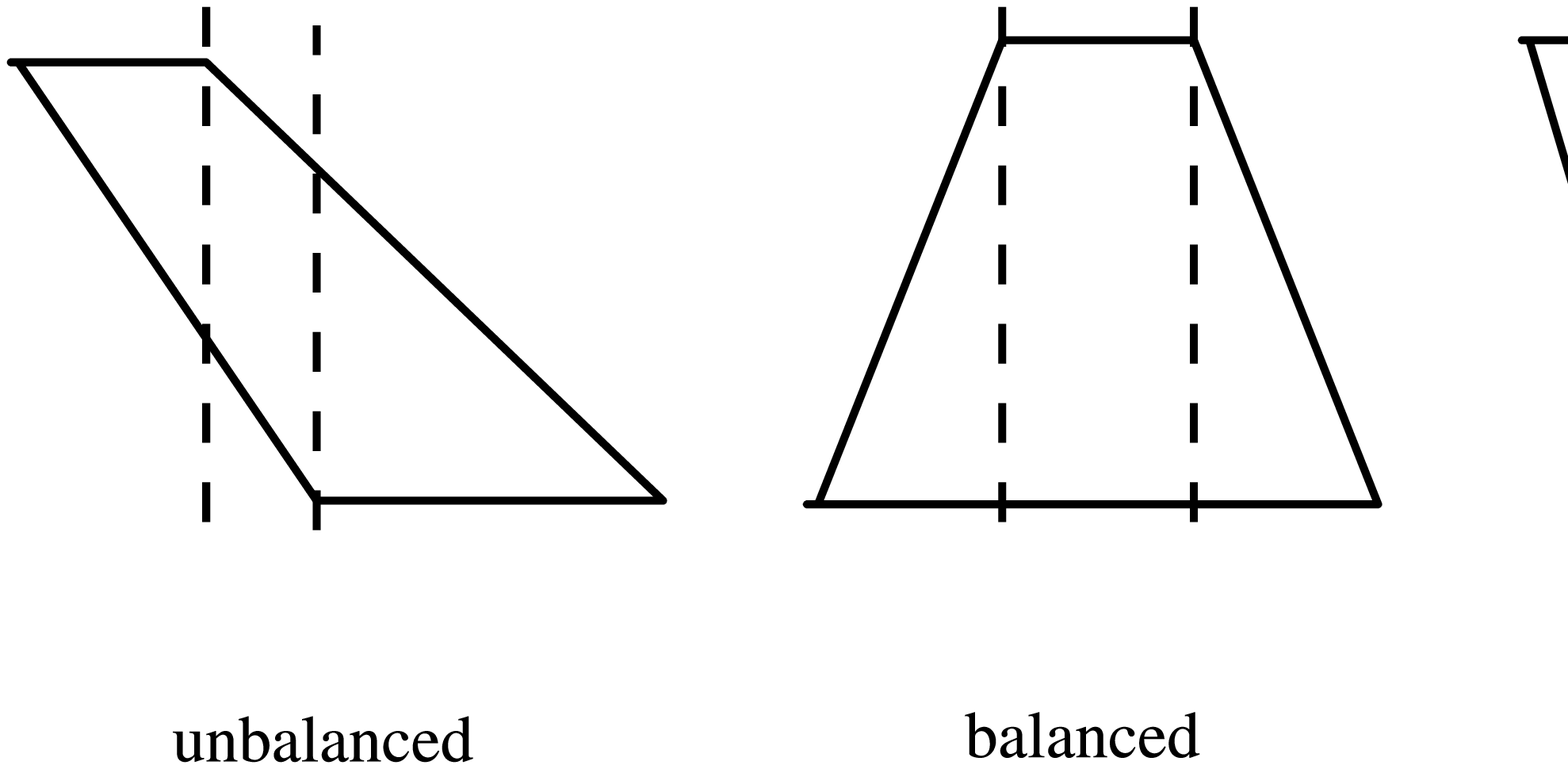}
       \caption{}
       \label{trapezes}
       \end{figure}

In the plane, we will call band of a line segment
\([\vec{a},\vec{b}]\) the set of the \(\vec{r}+\vec{u}\), where
\(\vec{r}\in [\vec{a},\vec{b}]\) and \(\vec{u}\bot \vec{b}-\vec{a}\).
We will call bands of a trapezoid the bands
associated with its bases. Let us denote them by \(\mathscr{B}\)
and \(\mathscr{B}'\).
If \(\mathscr{B} \subset \mathscr{B}'\)
or \(\mathscr{B}' \subset \mathscr{B}\), we will say that the trapezoid
is balanced. If \(\mathscr{B} \cap \mathscr{B}'=\emptyset\),
we will say that the trapezoid is unbalanced.
Otherwise, we will say that the trapezoid is semi-balanced (figure \ref{trapezes}).

\begin{thm}
\label{multnonnul}
For any two non-vanishing real numbers \(x\) and \(y\),
there exist exactly two non-collinear
central configurations (up to similarities) with non-vanishing multiplier
for the masses \(x\), \(-x\), \(y\), \(-y\).
These configurations are semi-balanced trapezoids.
They are symmetrical with respect to a line which is orthogonal to the bases.
If \(x=y\) or \(x=-y\), they are diamonds.
The bodies with masses \(x\) and \(-x\)
(respectively the bodies with masses \(y\) et \(-y\))
are at the endpoints of one of the two parallel sides.
Among the parallel sides, the side whose endpoints have the larger
absolute value is the smaller.
Two bodies which are on a same diagonal have masses with the same sign.
\end{thm}

{\bf Proof.} According to proposition \ref{inertiestat},
we have, for such a central configuration: \(\vec{\lambda}=\vec{0}\).
So \(x(\vec{r}_2-\vec{r}_1)=y(\vec{r}_3-\vec{r}_4)\):
the configuration is a trapezoid,
and the bodies with masses \(x\) and \(-x\)
(respectively the bodies with masses \(y\) et \(-y\))
are at the endpoints of one of the two parallel sides.
Among the parallel sides, the side whose endpoints have the larger
absolute value is the smaller.
Two bodies which are on a same diagonal have masses with the same sign.\\

We can assume: \(x\), \(y>0\).
If necessary, we can exchange body \(1\) and body \(2\),
or body \(3\) and body \(4\).\\

Let us assume the trapezoid to be balanced
and, for instance: \(y\le x\). Then we have:
\(||\vec{r}_4-\vec{r}_1||<||\vec{r}_4-\vec{r}_2||\),
\(||\vec{r}_3-\vec{r}_2||<||\vec{r}_3-\vec{r}_1||\).
This is incompatible with the identity:
\[\frac{1}{||\vec{r}_4-\vec{r}_2||^3}+\frac{1}{||\vec{r}_3-\vec{r}_1||^3}
=\frac{1}{||\vec{r}_4-\vec{r}_1||^3}+\frac{1}{||\vec{r}_3-\vec{r}_2||^3},\]
which follows from proposition \ref{cc4lambda0}.
So the trapezoid is not balanced.
With a similar reasoning, we prove that the trapezoid is not unbalanced.\\

As we are looking for configurations up to similarities, we can assume
\(||\vec{r}_4-\vec{r}_3||=x\). The configuration is a trapezoid.
Let us denote by \(\vec{a}\) the intersection of
\([\vec{r}_1, \vec{r}_3]\) and \([\vec{r}_2, \vec{r}_4]\),
which are the diagonals
of the trapezoid. The trapezoid up to isometries is defined by
\(||\vec{r}_3-\vec{a}||\) and \(||\vec{r}_4-\vec{a}||\)
such that \((||\vec{r}_3-\vec{a}||,||\vec{r}_4-\vec{a}||,x)\)
satisfies the triangular inequalities. Let us set:
\(u=||\vec{r}_3-\vec{r}_1||^2\), \(v=||\vec{r}_4-\vec{r}_2||^2\).
We have:
\(\frac{\sqrt{v}}{||\vec{r}_4-\vec{a}||}
=\frac{\sqrt{u}}{||\vec{r}_3-\vec{a}||}=\frac{x+y}{x}\).
The trapezoid up to isometries is defined by
\((u,v)\) such that \((\sqrt{u},\sqrt{v},x+y)\) are the lengths of the sides of a triangle.\\

Let us set:
\[I(\vec{a})=\sum _{i=1} ^4 m_i||\vec{r}_i-\vec{a}||^2\cdot \]
We have: \(dI=-2\vec{\lambda}.d\vec{a}=0\). Hence:
\[\left \{ \begin{array}{c}
I(\vec{r}_2)=I(\vec{r}_1)\\
I(\vec{r}_2)=I(\vec{r}_3) \end{array} \right .\]
This is equivalent to:
\[\left \{ \begin{array}{c}
||\vec{r}_3-\vec{r}_2||^2+||\vec{r}_4-\vec{r}_1||^2
=u+v-2xy\\
(x+y)||\vec{r}_3-\vec{r}_2||^2=xu+yv-xy(x+y)
\end{array} \right .\]
We set:
\[\varphi \left ( \begin{array}{c}
u\\
v\end{array} \right )
=\left ( \begin{array}{cc}
\frac{x}{x+y} & \frac{y}{x+y}\\
\frac{y}{x+y} & \frac{x}{x+y}
\end{array} \right )
\left ( \begin{array}{c}
u\\
v\end{array} \right )
-xy\left ( \begin{array}{c}
1\\
1\end{array} \right )\cdot \]
We have: \(\varphi =\mathscr{T} \circ \mathscr{A}\),
where \(\mathscr{A}\) is the
orthogonal scaling with respect to the first bisector
with scale factor \(\frac{x-y}{x+y}\) (it is an orthogonal
projection for \(x=y\)), whose absolute value is smaller than \(1\),
and \(\mathscr{T}\) is the translation with vector \(-xy(1,1)\).
The previous system is equivalent to:
\((||\vec{r}_3-\vec{r}_2||^2,||\vec{r}_4-\vec{r}_1||^2)
=\varphi (u,v)\).\\

\begin{figure}[h]
  \centering
     \includegraphics[width=7cm]{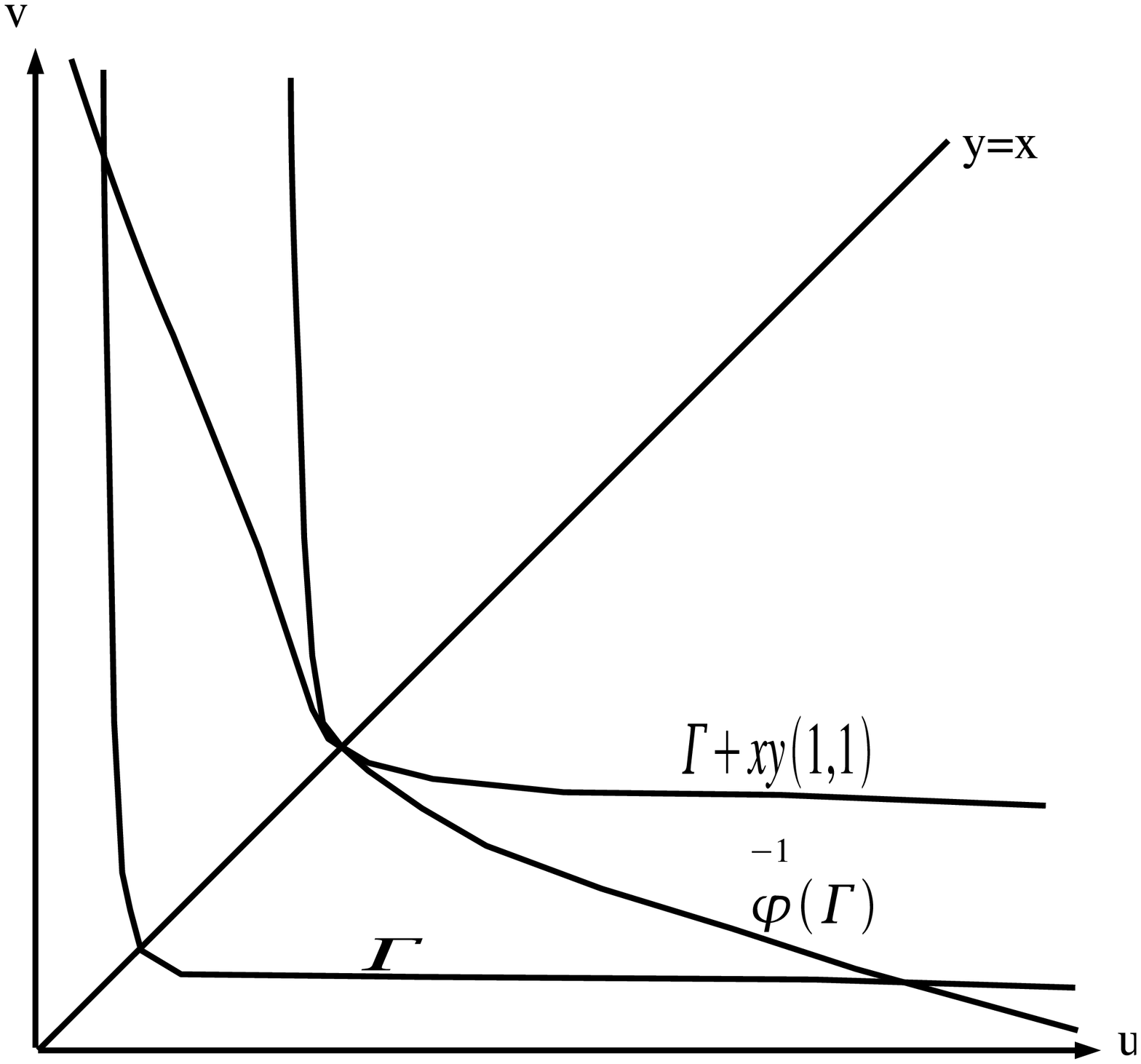}
       \caption{}
       \label{graphmultnon0}
       \end{figure}

If the configuration is non-collinear central, we have,
according to proposition \ref{cc4lambda0}:
\[\frac{1}{||\vec{r}_3-\vec{r}_2||^3}+\frac{1}{||\vec{r}_4-\vec{r}_1||^3}
=\frac{1}{u^{\frac{3}{2}}}+\frac{1}{v^{\frac{3}{2}}}
=\frac{1}{x^3}+\frac{1}{y^3}\cdot \]
We set: \(f(u,v)=\frac{1}{u^{3/2}}+\frac{1}{v^{3/2}}\).
The previous identity is equivalent to:
$$f(u,v)=f(\varphi (u,v))=f(x^2,y^2) \eqno(1)$$ 
Let \(\Gamma \) be the contour line of \(f\)
associated with the value
\(f(x^2,y^2)\). Then \((u,v)\) satisfies equation \((1)\)
if, and only if: \((u,v)\in \Gamma \cap \varphi ^{-1}(\Gamma)\).
For every \(x\), \(y>0\), there exist exactly two solutions of \((1)\)
in \(\mathbb{R}_+^*\times \mathbb{R}_+^*\), which can be written as \((u,v)\)
and \((v,u)\). The solutions are the two intersections of
\(\Gamma\) and \(\varphi ^{-1}(\Gamma)\), and they depend on \(x\) and
\(y\) in a continuous way (figure \ref{graphmultnon0}).\\

Conversely, let us consider one of the two solutions \((u,v)\)
of \((1)\). We have just seen that 
\((u,v)\) defines a trapezoid up to isometries
if, and only if, the numbers \((\sqrt{u},\sqrt{v},x+y)\) are the lengths 
of the sides of a triangle. If \(x=y=1\), we obtain:
\[(u,v)\approx (3,332979836;0,6670201635)
\mbox{ or } (0,6670201635;3,332979836)\cdot \]
Thus for any \(x=y\):
\[(\sqrt{u},\sqrt{v},x+y)\approx (1,825645047;0,8167130240;2)x\]
\[\mbox{ or } (0,8167130240;1,825645047;2)x\cdot \]
It can be checked that they are the sides of a non-flat triangle.
If \(x\ne y\),
let us assume that \((\sqrt{u},\sqrt{v},x+y)\) are not the lengths 
of the sides of a non-flat triangle. Let us denote by \((u(t),v(t))\)
a solution (depending on \(t\)) associated with
\(x(t)=x\) and \(y(t)=x+t(y-x)\).
We have \(x(0)=y(0)\), and
\((\sqrt{u(0)},\sqrt{v(0)},x(0)+y(0))\) actually defines a non-flat triangle,
whereas \((\sqrt{u(1)},\sqrt{v(1)},x(1)+y(1))\)
does not define a non-flat triangle any more.
The solutions \(u\) and \(v\) depend on \(x\) and
\(y\) in a continuous way. So at a certain time \(t\le 1\), the triangle is flat.
So the configuration of the four bodies at time \(t\) is a flat trapezoid.
As the trapezoid is semi-balanced,
the bodies are on a line in the order: \(1\), \(4\), \(2\),
\(3\). Hence:
\(||\vec{r}_4-\vec{r}_1||<||\vec{r}_2-\vec{r}_1||\) and
\(||\vec{r}_3-\vec{r}_2||<||\vec{r}_4-\vec{r}_3||\). Thus:
\[\frac{1}{||\vec{r}_2-\vec{r}_1||}+\frac{1}{||\vec{r}_4-\vec{r}_3||}
<\frac{1}{||\vec{r}_4-\vec{r}_1||}+\frac{1}{||\vec{r}_3-\vec{r}_2||},\]
which is impossible.
So \((\sqrt{u},\sqrt{v},x+y)=(\sqrt{u(1)},\sqrt{v(1)},x(1)+y(1))\)
actually defines a non-flat triangle.
So, according to proposition \ref{cc4lambda0}, the solution \((u,v)\)
defines a non-collinear central configuration
with vanishing vector of inertia. If \(x=y\), according to theorem
\ref{multzero}, we cannot have: \(\xi =0\).
According to the theorem, this is still true
if \(x\ne y\) as the bodies with masses \(x\) and \(-x\) cannot be at the endpoints
of a same diagonal.
Thus the solution \((u,v)\) actually defines a non-collinear
central configuration with vanishing multiplier.\\

We obtain a central configuration from the other by exchanging
\(u\) et \(v\). This is equivalent to exchanging body \(1\) and body \(2\)
and to exchanging body \(3\) and body \(4\).
This is also equivalent to considering the image of the configuration
by a symmetry whose axis is orthogonal to the bases of the trapezoid.
For one solution, the endpoints of the larger diagonal of the trapezoid
are the bodies with positive masses. For the other one, they are the bodies with negative masses.\\

For \(x=y\), every central configuration
is a parallelogram according to the vanishing of
\(\vec{\lambda}\). In fact, in this case, vector
\(\varphi (\tilde{u},\tilde{v})\) is collinear with \((1,1)\)
for every \((\tilde{u},\tilde{v})\), so
\((||\vec{r}_3-\vec{r}_2||^2,||\vec{r}_4-\vec{r}_1||^2)=\varphi (u,v)\)
is collinear with \((1,1)\). Hence:
\[\frac{2}{||\vec{r}_3-\vec{r}_2||^3}
=\frac{1}{||\vec{r}_3-\vec{r}_2||^3}+\frac{1}{||\vec{r}_4-\vec{r}_1||^3}
=\frac{1}{||\vec{r}_2-\vec{r}_1||^3}+\frac{1}{||\vec{r}_4-\vec{r}_3||^3}
=\frac{2}{x^3} \cdot \]
The configuration is a diamond. QED

\section{A result on co-circular configurations.}

\begin{thm}
\label{circulaire}
There is no co-circular four-body central configuration
with vanishing total mass and vector of inertia.
\end{thm}

{\bf Proof.} Let us consider a possible co-circular
four-body central configuration with vanishing vector of inertia.
We can assume that the bodies are in the order: 
\(1\), \(2\), \(3\), \(4\) on the circle.
We denote by \(\theta _{ijk}\) the length of the non-oriented circular arc
with ends \(\vec{r}_i\) and \(\vec{r}_k\) which contains
\(\vec{r}_j\). It belongs to \(]0, 2\pi[\).
As \(\theta _{123}+\theta _{143}=2\pi\), we can assume
\(\theta _{123} \le \pi\) (if necessary, we can exchange
body \(2\) and body \(4\)).
As \(\theta _{412}+\theta _{432}=2\pi\), we can assume
\(\theta _{432} \le \pi\) (if necessary, we can exchange
body \(1\) and body \(3\)).
As \(\theta _{123} \le \pi\), we have:
\(||\vec{r}_2-\vec{r}_1||<||\vec{r}_3-\vec{r}_1||\).
As \(\theta _{432} \le \pi\), we have:
\(||\vec{r}_4-\vec{r}_3||<||\vec{r}_4-\vec{r}_2||\).
Now proposition \ref{cc4lambda0} provides:
\[\frac{1}{||\vec{r}_2-\vec{r}_1||^3}+\frac{1}{||\vec{r}_4-\vec{r}_3||^3}
=\frac{1}{||\vec{r}_3-\vec{r}_1||^3}+\frac{1}{||\vec{r}_4-\vec{r}_2||^3},\]
which is impossible. QED\\

The study of the co-circular four-body central configurations
with non-vanishing total mass is difficult.
It is known that if the center of
inertia is the center of the circle, then the masses must be equal, and the configuration is a square ([7]).
We have the same result when the center of the circle is the intersection
of the diagonals ([12]).\\

The {\sl moment of inertia} of the configuration
with respect to the point \(\vec{a}\), which is denoted by \(I(\vec{a})\),
is defined by:
\[I(\vec{a})=\sum _{i=1} ^4 m_i||\vec{r}_i-\vec{a}||^2\cdot \]
In the proof of theorem \ref{multnonnul}, we saw that it was independent of \(\vec{a}\)
for \(\vec{\lambda}=\vec{0}\).
Thanks to Leibniz's formula, theorem \ref{circulaire} enables to show that
there is no non-collinear four-body central configuration
with vanishing total mass, vector of inertia and moment of inertia.

\section*{Acknowledgements.}
The results on configurations with \(\xi \ne 0\) are part of my Ph.D. thesis
([6]). I would like to express my gratitude to
my advisers Alain Chenciner and Alain Albouy
at the Institut de M\'ecanique C\'eleste of the Paris Observatory
for encouraging me to
study central configurations and systems with vanishing total mass,
for their advice, and for an elegant topological
argument in the proof of proposition \ref{dipole}.

\section*{References.}

[1] A. Albouy,
Sym\'etrie des configurations centrales de quatre corps.
C. R. Acad. Sci. Paris 320, pp. 217-220 (1995).\\

[2] A. Albouy,
The symmetric central configurations of four equal masses.
Contemp. Math. 198, pp. 131-135 (1996).\\

[3] A. Albouy,
On a paper of Moeckel on central configurations.
Regular Chaotic Dyn. 8, no 2, pp. 133-142 (2003).\\

[4] F. Alfaro, E. P\'erez-Chavela,
Families of continua of central configurations in charged problems.
Dyn. Contin. Discrete Impuls. Syst. series A: Mathematical Analysis 9,
pp. 463-475 (2002).\\

[5] M. Celli,
Homographic three-body motions with positive and negative masses.
In: Gaeta, G., Prinari, B., Rauch-Wojciechowski, S., Terracini, S. (eds.),
Symmetry and Perturbation Theory, pp. 75-82.
Proceedings of the International Conference
on SPT 2004, Cala Gonone, Italy, June 2004. World Scientific (2005).\\

[6] M. Celli,
Sur les mouvements homographiques de \(N\) corps associ\'es \`a des masses de signe quelconque, le cas particulier
o\`u la somme des masses est nulle, et une application \`a la recherche de
chor\'egraphies perverses. Ph.D. thesis, Paris 7 University (2005).
Can be found on the following website:\\
http://tel.ccsd.cnrs.fr/tel-00011790\\

[7] M. Hampton,
Co-circular configurations in the four-body problem.
In: Dumortier, F., Broer, H., Mawhin, J., Vanderbauwhede, A.,
Verduyn Lunel, S. (eds.),
Equadiff 2003, pp. 993-998.
Proceedings of the International Conference
on Differential Equations, Hasselt, Belgium, July 2003. World Scientific (2005).\\

[8] M. Hampton, R. Moeckel,
Finiteness of relative equilibria
of the four-body problem. Invent. Math. 163, no 2., pp. 289-312 (2006).\\

[9] E. Leandro,
Bifurcations and stability of some symmetrical
classes of central configurations. Ph.D. thesis, University of Minnesota (2001).\\

[10] M. Lindow,
Der kreisfall im problem der 3+1 k\"orper.
Astron. Nach. 220, pp. 369-380 (1924).\\

[11] G. Roberts,
A continuum of relative equilibria in the
five-body problem. Phys. D127, no 3-4, pp. 141-145 (1999).\\

[12] D. Saari,
Collisions, rings and other Newtonian \(N\)-body problems, p. 124.
Conference Board of the Mathematical Sciences, Regional Conference Series in Mathematics,
no 104. American Mathematical Society (2005).\\

[13] C. Sim\'o,
Relative equilibria in the four-body problem.
Cel. Mech. 18, pp. 165-184 (1978).\\

[14] S. Smale,
Mathematical problems for the next century.
Math. Intelligencer 20, no 2, pp. 7-15 (1998).

\end{document}